\documentclass[11pt]{amsart}

\usepackage{enumerate}

\numberwithin{equation}{section}

\newcommand{\bE}{\mathbb{E}}
\newcommand{\bK}{\mathbb{K}}
\newcommand{\bN}{\mathbb{N}}
\newcommand{\bP}{\mathbb{P}}
\newcommand{\bQ}{\mathbb{Q}}
\newcommand{\bR}{\mathbb{R}}
\newcommand{\bT}{\mathbb{T}}
\newcommand{\cP}{\mathcal{P}}

\newcommand{\ee}{\mathbf{e}}

\newcommand{\YY}{\mathbf{Y}}
\newcommand{\ZZ}{\mathbf{Z}}

\newcommand{\yy}{\mathbf{y}}
\newcommand{\zz}{\mathbf{z}}
\newcommand{\uu}{\mathbf{u}}
\newcommand{\vv}{\mathbf{v}}

\newcommand{\la}{\lambda}
\newcommand{\ep}{\epsilon}

\newcommand{\goto}{\rightarrow}

\theoremstyle{plain} \newtheorem{Theo}{Theorem}[section]
\theoremstyle{plain} \newtheorem{Lemma}[Theo]{Lemma}
\theoremstyle{plain} \newtheorem{Cor}[Theo]{Corollary}
\theoremstyle{plain} \newtheorem{Prop}[Theo]{Proposition}
\theoremstyle{remark} \newtheorem{Rem}[Theo]{Remark}
\theoremstyle{definition} 
\theoremstyle{definition} 
\theoremstyle{remark}

\begin{document}

\centerline{\large \bf Stepping-stone model with circular Brownian
migration}

\bigskip

\centerline{Xiaowen Zhou}

\bigskip

\centerline{Department of Mathematica and Statistics, Concordia
University}

\medskip

\centerline{Montreal, Canada}

\bigskip

\centerline{Abstract}
\medskip

In this paper we consider a stepping-stone model on a circle with
circular Brownian migration. We first point out a connection between
Arratia flow and the marginal distribution of this model. We then
give a new representation for the stepping-stone model using Arratia
flow and circular coalescing Brownian motion. Such a representation
enables us to carry out some explicit computation. In particular, we
find the Laplace transform for the time when there is only a single
type left across the circle.

\medskip



\section{Introduction}

Stepping-stone model is a mathematical model for population
genetics. A discrete-site stepping-stone model describes  the
simultaneous evolution of interacting populations over a collection
of finite or countable colonies. There are mutation, selection and
resampling within each colony, and there is migration among
different colonies. See \cite{Kim53} and \cite{Shi88} for some early
work.

Continuous-site stepping-stone model was first introduced in
\cite{Eva97}.  It is a process takeing values from the space
\[\Xi:= \{\mu: E \goto M_1(\bK)\},\] where $E$ denotes the continuous {\it site space},
$\bK$ denotes the {\it type space}, and $M_1(\bK)$ denotes the space
of all probability measures on $\bK$. Intuitively, such a map $\mu$
simultaneously represents the relative frequencies of different
types in populations at various sites. More precisely, for $e\in E$
and $B\subset \bK$, $\mu(e)(B)$ represents the ``proportion of the
population at the site $e$ processing types from the set $B$''. The
``moments'' of the continuous-site stepping-stone model are
specified using the so called {\it migration processes} taking
values in $E$.

In this paper we only consider a stepping-stone model with site
space $\bT$, a circle of circumference $1$, with type space
$\bK=[0,1]$, and with Brownian migration on $\bT$. We call it  a
stepping-stone model with circular Brownian migration (in short, a
SSCBM)  and write it as $X$ throughout the paper.

The distribution of SSCBM is uniquely determined by a family of
coalescing Brownian motions on $\bT$. But we have to go through more
notations before we could present the explicit formula.

Given a positive integer $n$, let $\cP_n$ denote the set of {\em
partitions} of $\bN_n := \{1, \ldots, n\}$.  That is, an element
$\pi$ of $\cP_n$ is a collection $\pi = \{A_1, \ldots, A_h\}$ of
disjoint subsets of $\bN_n$ such that $\bigcup_i A_i = \bN_n$. The
sets $A_1, \ldots A_h$ are the {\em blocks} of the partition $\pi$.
The integer $h$ is called the {\em length} of $\pi$ and is denoted
by $|\pi|$. Equivalently, we can think of $\cP_n$ as the set of
equivalence relations on $\bN_n$ and write $i \sim_\pi j$ if $i$ and
$j$ belong to the same block of $\pi \in \cP_n$.

Given $\pi\in \cP_n$, let
$$\alpha_i := \min A_i,\,  1\leq i\leq |\pi|.$$
$\{\alpha_i\}$ is the collection of {\it minimal elements} for
$\pi$.

By a circular (instantaneously) coalescing Brownian motion we mean a
collection of Brownian motions on $\bT$ such that any two of them
will move together as soon as they first meet. Given a circular
coalescing Brownian motion $(Z_1,\ldots,Z_n)$ starting at
$\ee=(e_1,\ldots,e_n)$. For $t>0$, let $\pi^\ee(t)$ be a
$\cP_n$-valued random partition such that $i\sim_{\pi^\ee(t)} j$ iff
$Z_i(t)=Z_j(t)$. Then $\pi^\ee(t)$ is the random partition {\it
induced} by $(Z_i(t))$. Write
\[\Gamma^\ee(t):=\{\alpha_i(t):1\leq i\leq |\pi^\ee(t)|\}\]
for the collections of  minimal elements for $\pi^\ee(t)$.

SSCBM is  then a $\Xi$-valued Hunt process $X$ with its transition
semigroup specified as following. Given $\mu\in\Xi$ and $n>0$, for
any $f_i\in C(\bT)$, $K_i\subset\bK$, $i=1,\ldots,n$,
\begin{equation}
\label{mom-dual}
\begin{split}
& \bQ^\mu\left[\prod_{i=1}^n \int_\bT de_i f(e_i)X_t(e_i)(K_i)
\right]
 = \int_{\bT^n}d\ee\prod_{i=1}^n f(e_i)\bP\left[
\bigotimes_{i\in \Gamma^\ee(t)}
\mu(Z_i(t))(\bigcap_{j\sim_{\pi^\ee(t)}i}K_j) \right],
\end{split}
\end{equation}
where $\bQ^\mu$ denotes the probability law of $X$ when its initial
value is $\mu$. See Theorem 4.1 in \cite{DEFKZ00} for a result on a
general continuous-site stepping-stone model.

In \cite{DEFKZ00} a particle representation for $X$ was given  using
the Poisson random measure on $D_\bT [0,\infty[\times \bK $ and a
``look down" scheme similar to that in \cite{DoKu96}. It leads to
better insight into the model. In the same spirit we are going to
propose another representation for $X$ in this paper.

It was shown in \cite{Eva97} that $X$ degenerates, i.e. for any
$t>0$, for almost all $e\in\bR$, $X_t(e)$ becomes a point mass on
some $k\in\bK$. A stronger version of this clustering behavior was
later shown in \cite{DEFKZ00} for site-space $\bT$ and in
\cite{Zho02} for site-space $\bR$. In fact, when the site-space is
$\bR$ there exists a random partition of $\bR$ such that $\bR$ is
divided into disjoint intervals and $X_t(e)$ is a point mass on the
same $k\in \bK$ for almost all $e$ in each interval. This suggests
that we can identify $X_t$ with a function $f$ on $\bT$ such that
$f(e)=k$ whenever $X_t(e)=\delta_{k}$. In this way we can identify
$X$ with  a step-function-valued process.

Let $\Xi'$ be the space of $\bK$-valued right continuous step
functions on $\bT$ equipped with the topology inherited from
$D_\bR(\bT)$. For each $\mu\in\Xi$, we are going to construct a
$\Xi'$-valued process $(X'_t, t>0)$ which can be regarded as SSCBM
with initial value $\mu$ under the above-mentioned identification.
To this end, we first point out an interesting connection between
 Arratia flow and SSCBM in Section 2. This connection allows us to specify
 the {\it entrance law} of $X'$ using the pre-image of Arratia flow. Then we  give an explicit
construction of $X'$ using Arratia flow and circular coalescing
Brownian motion, and we will show that $X'$ so defined does have the
right distribution under the above-mentioned identification. In this
sense $X'$ provides a nice version for $X$. Such a representation
enables us to compute the distribution of the time when there is
only a single type of individuals left across $\bT$ in Section 3. It
also allows us to obtain a result on the type that survives
eventually.

\section{A representation of stepping-stone model with circular Brownian migration}

We adopt some conventions for the rest of this paper. We identify
$\bT$ with interval $[0,1)$. Whenever we write $(e_1,\ldots,e_m)\in
\bT^m$ it implies that $e_1,\ldots,e_m$ have been already arranged
in anti-clockwise order around $\bT$. Given $u,v \in \bT$, write
$[u,v[$ for an interval starting at $u$ and ending at $v$ in
anti-clockwise order. Write $v-u$ for the length of the interval
$[u,v[. $ For $\{k_i\}\subset \bK$ and $(e_1,\ldots,e_m)\in\bT^m$,
write $\sum_{i=1}^m k_i1\{[e_i,e_{i+1}[\}, \, e_{m+1}:=e_1 $, for a
right continuous step function on $\bT$.

Arratia flow was first introduced in \cite{Arr79}. Arratia flow on
$\bT$ describes the evolution of a stochastic system in which there
is one Brownian motion starting at each point in $\bT$. Two Brownian
motions coalesce once they meet.  Formally, the Arratia flow can be
defined as a collection $\{\phi(s,t,x):0\leq s\leq t, x\in \bT\}$ of
random variables such that

\begin{itemize}
\item the random map $(s,t,x) \mapsto \phi(s,t,x)$ is jointly
measurable,
\item for each $s$ and $x$, the map $t \mapsto
\phi(s,t,x)$, $t \ge s$, is continuous,
\item for each $s$ and $t$
with $s \le t$, the map $x \mapsto \phi(s,t,x)$ is non-decreasing
and right-continuous, \item for $s \le t \le u$, $\phi(t,u,\cdot)
\circ \phi(s,t,\cdot) = \phi(s,u,\cdot)$,
\item for $u>0$,
$(s,t,x) \mapsto \phi(s+u, t+u, x)$ has the same distribution as
$\phi$, \item for $(x_1,\ldots,x_m)\in \bT^m$ the process
$(\phi(0,t,x_1), \ldots, \phi(0,t,x_m))_{t \ge 0}$ has the same
distribution as a circular coalescing Brownian motion starting at
$(x_1, \ldots, x_m)$.
\end{itemize}

From the continuity of Brownian  sample paths we see that, for each
$t>0$, there exists a positive integer valued random variable $N(t)$
and two sequences of random variables $(V_i(t))\in \bT^{N(t)}$ and
$(U_i(t))\in \bT^{N(t)}$ such that
\[\phi(0,t,x)=V_i(t) \,\,\text{for}\,\, x\in [U_i(t),U_{i+1}(t)[ \, \,
\text{and} \,\, i=1,\ldots,N(t).\]
In fact, we can even show that
\[\bP[|N(t)|]=1+2\sum_{n=1}^\infty \exp\left\{-n^2\pi^2 t\right\},\]
where $|N(t)|$ denotes the cardinality for $N(t)$. See Corollary 9.3
in \cite{DEFKZ00}.

For any $\mu\in\Xi$, given $N(t)$, $(U_i(t),i=1,\ldots,N(t))$ and
$(V_i(t),i=1,\ldots,N(t))$, let $\{\kappa_i, i=1,\ldots,N(t)\}$ be a
collection of independent $\bK$-valued random variables such that
$\kappa_i$ follows the distribution $\mu(V_i(t))$. Define
\begin{equation}\label{arratia1}
X'_t(e)=\sum_{i=1}^{N(t)}\kappa_i 1\{[U_i(t),U_{i+1}(t)[\}(e), \,\,
e\in\bT.
\end{equation}

We first point out that  $X'(t)$, when identified as
\begin{equation}\label{iden}
\sum_{i=1}^{N(t)}\delta_{\kappa_i} 1\{[U_i(t),U_{i+1}(t)[\},
\end{equation}
is indeed  a version of $X_t$.

\begin{Prop}\label{entrance}
 For any $t>0$, with the identification (\ref{iden}) $X'_t$ has the same distribution as $X_t$ under
$\bQ^\mu$.
\end{Prop}

\begin{proof}
To determine the distribution of $X'_t$ we only need to specify
joint distributions such as
 \[\bP\{X'_t(e_1)\in
dk_1,\ldots,X'_t(e_n)\in dk_n\}.\]

By definition $(\phi(0,t,e_1),\ldots,\phi(0,t,e_n))$ is a circular
coalescing Brownian motion starting at $(e_1,\ldots,e_n)$. Let
$\pi^\ee(t)$ be the induced partition on $\bN_n$. For any $k_i\in
\bK, i=1,\ldots,n$, given $(U_i(t))$  and $(V_i(t))$ as before,
observe that $\phi(0,t,e_i)$ and $\phi(0,t,e_j)$ belong to the same
interval $[U_r(t), U_{r+1}(t)[$ for some $r$ iff $i\sim_{\pi^\ee(t)}
j$. Also notice that $\phi(0,t,e_i)=\phi(0,t,e_j) $ whenever
$i\sim_{\pi^\ee(t)} j$. Then
\begin{equation}
\label{mom-dual-circle}
\begin{split}
& \bP\left\{\bigcap_{i=1}^n  \{X'_t(e_i)\in dk_i\} \right\} =
\bP\left[ \prod_{i\in
\Gamma^\ee(t)}1\left\{\bigcap_{j\sim_{\pi^\ee(t)}i} \{k_j=k_{i}\}
\right\} \mu\left(\phi(0,t,e_{i})\right)\left(dk_{i}\right) \right].
\end{split}
\end{equation}
An inspection of (\ref{mom-dual-circle}) reveals that
(\ref{mom-dual}) holds for $X'(t)$ when it is regarded as
$\Xi$-valued. So $X_t$ and $X'_t$ have the same distribution.

\end{proof}

By Proposition \ref{entrance} we may and will suppose that $X_t,
t>0$, is $\Xi'$-valued in the rest of the paper.

 We can read off some properties for $X_t$, $t>0$, immediately from
Proposition \ref{entrance}. First, with probability one $X_t$ (as a
function of $e$) can only take finitely many different values from
$\bK$. Moreover, if $\mu(e)$ is a diffuse measure for almost all
$e\in\bT$, then with probability one $X_t$  takes different values
over different intervals on $\bT$, i.e. $X(e)=k$ for $e\in
[e_1,e_2[$ whenever $X(e_1)=k=X(e_2)$. Such properties are also
discussed in Section 10 of \cite{DEFKZ00}.

Conditioning on $X_s=\sum_{i=1}^m k_i 1\{[u_i, u_{i+1}[\}$,
(\ref{mom-dual}) shows that, given $t>s$,  $X_t$ can only take
values from $\{k_i\}$. Moreover, for any $\{k'_j,
j=1,\ldots,n\}\subset \{k_i\}$ and any $(z_j)\in \bT^n$, by
(\ref{mom-dual}) we can further show that
\begin{equation}
\label{mom-dual2}
\begin{split}
& \bQ\left\{\left.\bigcap_{j=1}^n \{ X_t(z_j)=k'_j \}\right|
X_s=\sum_{i=1}^m k_i 1\{[u_i, u_{i+1}[\} \right\}\\
& \quad =\bQ\left\{\bigcap_{j=1}^n \{ X_s(Z_j(t-s))=k'_j \}\right\},
\end{split}
\end{equation}
where $(Z_j)$ is a circular coalescing Brownian motion starting at
$(z_j)$.

To describe the evolution of $X$ over time we need a Lemma on
duality between two circular coalescing Brownian motions.

Fix $\yy=(y_1,\ldots,y_m)\in \bT^m$  and $\zz=(z_1,\ldots,z_n)\in
\bT^n$. Let $(Y_1,\ldots,Y_m)$ be an $m$-dimensional circular
coalescing Brownian motion starting at $\yy$. Let $(Z_1,\ldots,Z_n)$
be an $n$-dimensional circular coalescing Brownian motion starting
at $\zz$. Put
\[
I_{ij}^\rightarrow(t,\zz) := 1\{Y_i(t) \in [z_{j}, z_{j+1}[\}
\]
and
\[
I_{ij}^\leftarrow(t,\yy) := 1\{y_i \in [Z_{j}(t), Z_{j+1}(t)[\}
\]
for $1 \le i \le m$ and $1\le j \le n$. Recall that $z_{n+1}:=z_1$
and $Z_{n+1}:=Z_1$.

\begin{Lemma}\label{coal-dual}
 The two $(m\times n)$-dimensional  arrays $(I_{ij}^\leftarrow(t,\yy))$
and $(I_{ij}^\rightarrow(t,\zz))$ have the same distribution.
\end{Lemma}

\begin{proof}
We can prove Lemma \ref{coal-dual} in the same way as Theorem 2.1 in
\cite{Zho05}, i.e. we first show that the corresponding duality
holds for circular coalescing random walks, and then apply
time-space scaling to obtain the desired result for circular
coalescing Brownian motions.

\end{proof}

By Lemma \ref{coal-dual} we can easily derive  the following side
result concerning a dual relationship for Arratia flow. Such a
result was pointed out in \cite{Arr79} for coalescing Brownian flow
on the real line.

\begin{Prop}
When identified as point processes, $(U_i(t))$ and $(V_i(t))$ have
the same distribution for any fixed $t>0$.
\end{Prop}

\begin{proof}
For any $(z_i)\in \bT^{2n}$, let $(Z_i)$ be a coalescing Brownian
motion starting at $(z_i)$. Consider a sequence of circular
coalescing Brownian motions $\{(\phi(0,t,x_i^m))_{i=1}^m,
m=1,2,\ldots\}$ such that the set $\{x_i^m, i=1,\ldots,m\}$ of
starting locations approaches to a dense set in $\bT$ as $m\goto
\infty$. Such a sequence provides an ``approximation'' for the
Arratia flow on $\bT$. Then by Lemma \ref{coal-dual}
\begin{equation}\label{arratia-dual}
\begin{split}
&\bP\left\{\bigcap_{i=1}^m\left\{\phi(0,t,x_i^m)\not\in
\bigcup_{j=1}^n (z_{2j-1},z_{2j})\right\}\right\}\\
&\quad=\bP\left\{\bigcap_{i=1}^m\left\{x_i^m\not\in \bigcup_{j=1}^n
(Z_{2j-1}(t),Z_{2j}(t))\right\}\right\}.
\end{split}
\end{equation}

 On the one hand, taking  limits on both sides of (\ref{arratia-dual}) as $m\goto\infty$, we can
show that
\[\bP\left\{\{V_i(t)\}\bigcap\bigcup_{j=1}^n (z_{2j-1},z_{2j})=\emptyset\right\}=
\bP\left\{\bigcap_{j=1}^n \left\{Z_{2j-1}(t)=Z_{2j}(t)
\right\}\right\}.\] On the other hand,
$\{U_i(t)\}\cap(z_{2j-1},z_{2j})=\emptyset$ iff
$(z_{2j-1},z_{2j})\subset [U_i(t),U_{i+1}(t)[$ for some $i$ iff
$\phi(0,t,z_{2j-1})=V_i(t)= \phi(0,t,z_{2j})$ for some $i$.
Consequently, we also have
\[\bP\left\{\{U_i(t)\}\bigcap\bigcup_{j=1}^n (z_{2j-1},z_{2j})=\emptyset\right\}=
\bP\left\{\bigcap_{j=1}^n \{\phi(0,t,z_{2j-1})=\phi(0,t,z_{2j})\}
\right\}.\]

Therefore, $(U_i(t))$ and $(V_i(t))$ have the same {\it avoidance
function}. So, the assertion of this Proposition holds (see Theorem
3.3 of \cite{Kla76}).
\end{proof}

Let us go back to the stepping-stone model. We first consider a
special initial value $\mu$. Given
\[\nu=\sum_{i=1}^m \delta_{k_i} 1\{[u_i,u_{i+1}[\}\in\Xi,\] write
$\YY=(Y_i)$ for an $m$-dimensional circular coalescing Brownian
motion starting at $\uu:=(u_i)$ and define
\begin{equation*}
X'_t=\sum_{i=1}^m k_i 1\{[Y_i(t), Y_{i+1}(t)[\}, \, t\geq 0,
\end{equation*}
with the convention that $1\{[y,y[\}:=0$.

\begin{Lemma}\label{repre1}
$X'$ has the same distribution as $X$ under $\bQ^\nu$.
\end{Lemma}

\begin{proof}

$(X'_t)$ is clearly a Markov process from its definition.

Given $\{k'_1,\ldots,k'_n\}\subset \{k_i, i=1,\ldots,m\}$ and
$(v_j)\in \bT^n$, let $\ZZ=(Z_j)$ be a circular coalescing Brownian
motion starting at $\vv:=(v_j)$. Set
\[g(\yy;\zz):=\prod_{j=1}^n\sum_{i:k_i=k'_j}1\{[y_{i-1},y_i[\}(z_j), \,\, \yy:=(y_i), \zz:=(z_j).\]
Lemma \ref{coal-dual} yields that
\begin{equation*}\label{eq1}
\begin{split}
\bP\left\{\bigcap_{j=1}^n\{X'_t(v_j)=k'_j\}\right\} &=
\bP\left\{\bigcap_{j=1}^{n} \left\{\sum_{i=1}^m k_i
1\{[Y_i(t), Y_{i+1}(t)[\}(v_j)=k'_j\right\} \right\}\\
&=\bP\left[g(\YY(t);\vv)\right]\\
&=\bP\left[g(\uu;\ZZ(t))\right]\\
&= \bP\left\{ \bigcap_{j=1}^{n} \left\{\sum_{i=1}^m k_i
1\{[u_i, u_{i+1}[\}(Z_j(t))=k'_j\right\} \right\}\\
&=\bP\left\{\bigcap_{j=1}^n\left\{X'_0(Z_j(t))=k'_j\right\}\right\}.
\end{split}
\end{equation*}
It then follows from Proposition \ref{entrance} and
(\ref{mom-dual2}) that $X'$ and $X$ have both the same initial value
and the same transition semigroup. So, they have the same
distribution.

\end{proof}

Now we are ready to construct a representation for $X$ with a
general initial value $\mu\in \Xi$. Given $\ep>0$, as in
(\ref{arratia1}) put
\begin{equation}\label{arratia3}
X'_\epsilon=\sum_{i=1}^{N(\ep)} \kappa_i 1\{[U_i(\ep),
U_{i+1}(\ep)[\}.
\end{equation}
 Given $N(\ep)$, $(U_1(\ep),\ldots,U_{N(\ep)}(\ep))$
and $(\kappa_1,\ldots,\kappa_{N(\ep)})$, write $(Y_i)$ for an
$N(\ep)$-dimensional circular coalescing Brownian motion starting at
$(U_i(\ep))$. We further define
\begin{equation}\label{arratia2}
X'_t=\sum_{i=1}^{N(\ep)} \kappa_i 1\{[Y_i(t-\epsilon),
Y_{i+1}(t-\epsilon)[\}, \, t\geq \ep,
\end{equation}
again, with the convention that $1\{[y,y[\}:=0$. Combining
Proposition \ref{entrance} and Lemma \ref{repre1} we can easily
obtain the following result.

\begin{Theo}\label{repre2}
Given $\mu\in\Xi$ and $\ep>0$, $(X'_t,t\geq \ep)$ has the same
distribution as $(X_t,t\geq \ep)$ under $\bQ^\mu$.
\end{Theo}

\begin{Rem}
The representation (\ref{arratia2}) suggests that SSCBM can also be
thought of as  a multi-type, nearest-neighbored voter model on
$\bT$. See Chapter V in \cite{Lig85} for discussions on voter model.
\end{Rem}

\begin{Rem}
A similar representation can be found for a stepping-stone model
with Brownian migration on $\bR$. We leave the details to the
readers.
\end{Rem}

\section{The first time when there is only a single type left}

In  this section we are going to study  properties of $X$ using the
representation given in Section 2.

Treating $(X_t,t>0)$ as $\Xi'$-valued, put
$$T:=\inf\{t> 0:\exists k\in \bK, X_t(e)=k ,\forall e\in \bT\}.$$
$T$ is then the first time when a single type of individuals prevail
all over $\bT$. It is easy to see from the representation
(\ref{arratia2}) that
\begin{equation*}\label{first1}
\bQ^\mu\{T<\infty\}=1,
\end{equation*}
for all $\mu\in\Xi$. Now we are going to find the exact distribution
for $T$.

We start with a preliminary result which is interesting in its own
right. Let $(Y_i)$ be an $m$-dimensional circular coalescing
Brownian motion starting at $(y_i)\in\bT^m$, $m\geq 2$. Let
$$T_m:=\inf\{t>0: Y_1(t)=\ldots=Y_m(t)\}. $$

\begin{Prop}\label{circle1}
Given any positive integer $m\geq 2$, we have
\begin{equation}\label{circle}
\bP[e^{-\lambda T_m}]=\sum_{i=1}^m
\frac{\sinh((y_{i+1}-y_i)\sqrt{\la})}{\sinh(\sqrt{\la})}, \,\,
\la>0.
\end{equation}
\end{Prop}

\begin{proof}
For $i=1,\ldots,m$, write $S_i$ for the time when $Y_{i+1}$ first
reaches $Y_i$ from the clockwise direction. As usual, we define
$Y_{m+1}:=Y_1$. Since $(Y_{i+1}-Y_i)/\sqrt{2}$ is again a Brownian
motion which starts at $(y_{i+1}-y_i)/\sqrt{2}$ and stops whenever
it reaches $0$ or $1/\sqrt{2}$, $S_i$ is then the first time that
the Brownian motion $(Y_{i+1}-Y_i)/\sqrt{2}$ reaches $1/\sqrt{2}$
before it reaches $0$. We thus have
\begin{equation}\label{time}
\bP[e^{-\lambda S_i}]=
\frac{\sinh((y_{i+1}-y_i)\sqrt{\la})}{\sinh(\sqrt{\la})}.
\end{equation}
See Exercise II.3.10 in \cite{ReYo91}.

Our key observation is that
\[\bP\{T_m<t\}=\bigcup_{i=1}^m\{S_i<t\},\]
and the events on the right hand side of this equation are disjoint.
So, (\ref{circle}) follows.

\end{proof}

Standard argument gives the following result.

\begin{Cor}
For any positive integer $m\geq 2$, we have
\[\bP[T_m]=\frac{1}{4}-\frac{1}{4} \sum_{i=1}^m (y_{i+1}-y_i)^3.\]
Consequently, $\bP[T_m]$ attains its maximum ${1}/{4}-{1}/{4m^2} $
iff all the initial values $y_1,\ldots,y_m$ are equally spaced on
$\bT$.
\end{Cor}

\begin{Rem}
An explicit expression for the distribution of $T_m$ can also be
found. By Theorem 4.1.1 in \cite{Kni81}, we have
\begin{equation*}\label{coal1}
\begin{split}
\bP\{S_i\leq
t\}=\sqrt{2}(y_{i+1}-y_i)+\frac{2}{\pi}\sum_{n=1}^\infty
\frac{(-1)^n}{n}\sin\left(\sqrt{2}n\pi(y_{i+1}-y_i)\right)
\exp\left\{-n^2\pi^2 t\right\}.
\end{split}
\end{equation*}
Therefore,
\begin{equation}
\begin{split}
\bP\{T_m\leq t\}=\sqrt{2}+\frac{2}{\pi}\sum_{i=1}^m\sum_{n=1}^\infty
\frac{(-1)^n}{n}\sin\left(\sqrt{2}n\pi(y_{i+1}-y_i)\right)
\exp\left\{-n^2\pi^2 t\right\}.
\end{split}
\end{equation}

\end{Rem}

We expect that $ \bP\{T_m\leq t\}$ also reaches its minimum when
$y_1, \ldots, y_m$ are equally spaced on $\bT$. But we do not have a
proof yet.

Let $\kappa$ be the type of individuals left after time $T$. Then
\[\kappa=\lim_{t\goto \infty} X_t(e), \, \forall e\in\bT.\]

\begin{Theo}\label{one-type}
Given $\mu\in \Xi$ such that $\mu(x)$ is a diffuse probability
measure for almost all $x\in\bT$, then the Laplace transform for $T$
has the expression
\begin{equation}\label{onetype1}
\bQ^\mu [e^{-\la T}]=\frac{\sqrt{\la}}{\sinh(\sqrt{\la})}, \,\,
\la>0.
\end{equation}
Moreover,
\begin{equation}\label{onetype2}
\bQ^\mu\{\kappa \in dk\}=\int_\bT de\mu(e)(dk).
\end{equation}
\end{Theo}

\begin{proof}
Given $m$ and  $(y_i)\in \bT^m$, we first observe that, by
(\ref{circle}),
\begin{equation}\label{first2}
\lim_{m\goto \infty}\bP[e^{-\lambda T_m}]=\lim_{m\goto
\infty}\sum_{i=1}^m
\frac{(y_{i+1}-y_i)\sqrt{\la}}{\sinh(\sqrt{\la})}=\frac{\sqrt{\la}}{\sinh(\sqrt{\la})}
\end{equation}
as \[\max_{1\leq i\leq m}(y_{i+1}-y_i)\rightarrow 0+.\]

Recall from the representation (\ref{arratia2}) that, for $\ep>0$,
\[X'_t=\sum_{i=1}^{N(\ep)} \kappa_i 1\{[Z_i(t-\epsilon), Z_{i+1}(t-\epsilon)[\},\,\, t\geq\ep, \]
where, given $N(\ep)$, $(Z_i)$ is a circular coalescing Brownian
motion starting at $(U_i(\ep))\in \bT^{N(\ep)}$. Put
\[T(\ep):=\inf\{t\geq 0: Z_1(t)=\ldots=Z_{N(\ep)}(t)\}.\]
Notice that, given $N(\ep)$, $\kappa_1,\ldots,\kappa_{N(\ep)}$ are
all different since $\mu(e)$ is diffuse for almost all $e\in\bT$.
Then $\ep+T(\ep)$ is also the first time when $X'(e)$  assumes
 a single value in $\bK$ for all $e\in\bT$.

Put
\[\Delta(\ep):=\max_{1\leq i\leq N(\ep)}(U_{i+1}(\ep)-U_i(\ep)).\]
It is evident from the definition of Arratia flow and the
representation (\ref{arratia2}) that $\Delta(\ep)\goto 0$ in
probability and
\[\bP\{T(\ep)>0\}\goto 1 \text{\,\, as \,\,} \ep\goto 0+.\] In
addition, $\bQ^\mu\{T>0\}=1$. It follows from Theorem \ref{repre2}
and (\ref{first2}) that
\begin{equation}
\begin{split}
\bQ^\mu[e^{-\la T}]&=\lim_{\ep\goto 0+}\bQ^\mu [e^{-\la T}; T>\ep] \\
&=\lim_{\ep\goto 0+}\bP\left[\bP\left[\left.e^{-\la
(\ep+T(\epsilon))}; T(\ep)>0\right|
 X'_\ep\right]\right] \\
&=\frac{\sqrt{\la}}{\sinh(\sqrt{\la})}.
\end{split}
\end{equation}

Finally, by (\ref{mom-dual}) we have
\[\lim_{t\goto\infty}\bQ^\mu\{X_t(e)\in dk\}=\lim_{t\goto\infty}\bP [\mu{(Z(t))}(dk)],  \]
where $Z$ is a circular Brownian motion starting at $e\in\bT$.
(\ref{onetype2}) thus follows.

\end{proof}

\begin{Rem}
Notice that the distribution of $T$ does not depend on $\mu$ as long
as $\mu(x)$ is diffuse for almost all $x\in\bT$.
\end{Rem}

\begin{Rem}
Let $X$ be a Brownian motion starting at $0<x<1/\sqrt{2}$. Put
\[T_x:=\inf\left\{t\geq 0: X_t=0 \text{\,\,or\,\,}{1}/{\sqrt{2}}\right\}.\]
We observe that, for the $\mu$ in Theorem \ref{one-type},
\begin{equation}\label{conditional}
\begin{split}
\bQ^\mu[e^{-\la T}]&=\lim_{x\rightarrow
0+}\frac{\sinh(x\sqrt{2\la})}{x\sqrt{2}\sinh(\sqrt{\la})}\\
&=\lim_{x\rightarrow 0+}\frac{\bP\left[e^{-\la
T_x};T_x={1}/{\sqrt{2}}\right]}
{\bP\left\{T_x={1}/{\sqrt{2}}\right\}}\\
&=\lim_{x\rightarrow 0+}\bP\left[e^{-\la T_x}\left|X_{T_x}={1}/{\sqrt{2}}\right.\right].\\
\end{split}
\end{equation}

Using (\ref{conditional}) and Theorem 4.1.1 in \cite{Kni81}  we can
further find an explicit expression for $\bQ^\mu\{T\leq t\}$. For
$t>0$,
\begin{equation}\label{distribution}
\begin{split}
\bQ^\mu\{T\leq t\} &=\lim_{x\goto 0+} \frac{1}{\sqrt{2}x}
\left(\sqrt{2}x+\frac{2}{\pi}\sum_{n=1}^\infty \frac{(-1)^n}{n}\sin
(\sqrt{2}n\pi x)
\exp\left\{-n^2\pi^2 t\right\}\right)\\
&=1+ 2\sum_{n=1}^\infty (-1)^n \exp\left\{-n^2\pi^2t\right\}.
\end{split}
\end{equation}

\end{Rem}

\begin{Rem}
It is not hard to see from the proof for Theorem \ref{one-type} that
the distribution (\ref{distribution}) coincides with the
distribution of the time when the image of Arratia flow on $\bT$
first becomes a set of a single element, i.e. the distribution of
\[\tau:=\inf\{t\geq 0: \phi(0,t,x)=\phi(0,t,y), \forall x, y\in\bT \}.\]

\end{Rem}

Again, for the $\mu$ given in Theorem \ref{one-type}, for any
$\ep>0$, let interval $[U'_\ep, U''_\ep[$ be the unique interval
$[U_i(\ep),U_{i+1}(\ep)[$ in (\ref{arratia3}) such that
$\kappa=k_i$; i.e. $[U'_\ep, U''_\ep[$ is the collection of sites at
time $\ep$ whose type eventually prevails.

\begin{Prop}
For the $\mu$ given in Theorem \ref{one-type}, as $\ep\goto 0+$ both
$(U'_\ep)$ and $(V''_\ep)$ converge in distribution to a uniform
distribution on $\bT$.
\end{Prop}

\begin{proof}
Clearly $U''_\ep-U'_\ep\goto 0 $ in probability. Therefore, we just
need to show that for any $[a,b[\subset \bT$,
\begin{equation}\label{last}
\lim_{\ep\goto 0+}\bP\{U'_\ep \in [a,b[\}=b-a.
\end{equation}
 To prove (\ref{last}), we first notice
that, given $N(\ep)$, events $\{S_i<\infty\}, 1\leq i\leq N(\ep) $,
are all disjoint, where $S_i$ is defined as in the proof for
Proposition \ref{circle1}, but for a coalescing Brownian motion
starting at $(U_i(\ep),U_{i+1}(\ep))$. Consequently, by (\ref{time})
\begin{equation*}
\begin{split}
\bP\{U'_\ep\in [a,b[\}&=\bP\left\{\bigcup_{1\leq i\leq
N(\ep)}\{U_i(\ep)\in [a,b[, S_i<\infty\}\right\}\\
&= \bE\left\{\sum_{1\leq i\leq N(\ep)}(U_{i+1}(\ep)-U_i(\ep)) 1\{U_i(\ep)\in [a,b[\}\right\}\\
&=\bE\left[U_n(\ep)-U_m(\ep)\right],
\end{split}
\end{equation*}
where  $m:=\min\{i:U_i(\ep)\in [a,b[\}$ and $n:=\max\{i:U_i(\ep)\in
[a,b[\}$. Therefore, (\ref{last}) follows readily.

\end{proof}

\begin{Rem}
If  $\mu\in\Xi$  is arbitrary, we can not find the explicit
distribution for $T$ under $\bQ^\mu$ . Nevertheless, similar to the
proof for Theorem \ref{one-type} we can still show that
\[\bQ^\mu\{T\leq t\}\geq 1+ 2\sum_{n=1}^\infty (-1)^n \exp\left\{-n^2\pi^2t\right\}, \,\, t\geq 0.\]
\end{Rem}

\bibliographystyle{alpha}
\bibliography{circle}

\end{document}